\documentclass[12pt]{amsart} 
\usepackage{amssymb}
\usepackage{multirow}
\usepackage{mdwlist}
\usepackage{enumerate}

\newtheorem{thm}{Theorem}[section]
\newtheorem{defin}[thm]{Definition}
\newtheorem{cor}[thm]{Corollary}
\newtheorem{prop}[thm]{Proposition}
\newtheorem{lemma}[thm]{Lemma}
\newtheorem{eg}[thm]{Example}
\newtheorem{rmk}[thm]{Remark}

\newcommand{\pf}[1][]{\emph{Proof#1:\ }}

\newcommand{\quadricaction}[2]{\ensuremath{\left\langle#1\,|\,#2\right\rangle}}

\newcommand{\integers}{\ensuremath{\mathbb{Z}}}
\newcommand{\rationals}{\ensuremath{\mathbb{Q}}}
\newcommand{\reals}{\ensuremath{\mathbb{R}}}
\newcommand{\complex}{\ensuremath{\mathbb{C}}}

\newcommand{\abs}[1]{\ensuremath{\left|#1\right|}}
\newcommand{\defsetshort}[1]{\ensuremath{\left\{#1\right\}}}
\newcommand{\defset}[2]{\defsetshort{#1\,\left|\,#2\right.}}
\newcommand{\defsetspan}[1]{\ensuremath{<\!\!#1\!\!>}}
\newcommand{\defsetspanlong}[2]{\defsetspan{#1\,\left|\,#2\right.}}

\newcommand{\defby}{\mathrel{\mathop:}=}

\newcommand{\symgroup}[1]{\ensuremath{\mathbb{S}_{#1}}}
\newcommand{\altgroup}[1]{\ensuremath{\mathbb{A}_{#1}}}
\newcommand{\dihgroup}[1]{\ensuremath{\mathbb{D}_{#1}}}
\newcommand{\cyclgroup}[1]{\ensuremath{\integers_{#1}}}
\newcommand{\bsymgroup}[1]{\ensuremath{\mathcal{S}_{#1}}}
\newcommand{\baltgroup}[1]{\ensuremath{\mathcal{A}_{#1}}}
\newcommand{\bdihgroup}[1]{\ensuremath{\mathcal{D}_{#1}}}
\newcommand{\bcyclgroup}[1]{\ensuremath{\mathcal{Z}_{#1}}}
\newcommand{\glgroup}[2]{\ensuremath{\mbox{GL}_{#1}\left(#2\right)}}
\newcommand{\slgroup}[2]{\ensuremath{\mbox{SL}_{#1}\left(#2\right)}}
\newcommand{\pglgroup}[2]{\ensuremath{\mbox{PGL}_{#1}\left(#2\right)}}
\newcommand{\Aut}[2][]{\ensuremath{\mbox{Aut}_{#1}\left(#2\right)}}

\newcommand{\proj}[2][]{\ensuremath{\mathbb{P}_{#1}^{#2}}}

\newcommand{\deffun}[2]{\ensuremath{#1\rightarrow#2}}
\newcommand{\deffunname}[3]{\ensuremath{#1:#2\rightarrow#3}}
\newcommand{\defmap}[2]{\ensuremath{#1\mapsto#2}}
\newcommand{\defmapname}[3]{\ensuremath{#1:#2\mapsto#3}}

\newcommand{\canclass}[1]{\ensuremath{\mbox{K}_{#1}}}

\title{Weakly-exceptional quotient singularities}
\author[Dmitrijs Sakovics]{Dmitrijs Sakovics}
\begin{document}
\maketitle

\begin{abstract}
A singularity is said to be weakly-exceptional if it has a unique purely log
terminal blow up. In dimension $2$, V.~Shokurov proved that
weakly-exceptional quotient singularities are exactly those of types $D_{n}$,
$E_{6}$, $E_{7}$, $E_{8}$. This paper classifies the weakly-exceptional
quotient
singularities in dimensions $3$ and $4$.
\end{abstract}

\setcounter{tocdepth}{1}
\tableofcontents
\section{Introduction}

Let $G\subset\glgroup{N}{\complex}$ be a finite subgroup. Then $G$ has a natural
linear action on $\complex^{N}$. The aim of this paper is to study the
exceptionality properties of quotient singularities $\complex^{N}/G$ in low
dimensions.

The Chevalley-Shephard-Todd theorem (see~\cite[Theorem~4.2.5]{Springer77}) easily implies that
as far as the exceptionality of the induced quotient singularities is
concerned, one can assume that $G$ contains no quasi-reflections.
Furthermore, the exceptionality of the singularity is only dependent of the
group's induced action on \proj{N-1} (for details, see
Section~\ref{sect:criteria}), so one can assume that $G$ is actually a subgroup of
\slgroup{N}{\complex}. Let $\bar{G}$ be the image of $G$ in the natural
projection of \slgroup{N}{\complex} onto \pglgroup{N}{\complex}. Since the two
groups are closely related, the argument will move from considering $\bar{G}$ to
considering $G$ and back at will. 

Standard notation will be used throughout this paper, but for the convenience
of the reader and in order to prevent any ambiguity, the notation used will be
fully described in Section~\ref{notation}.

Since the singularities are defined by group actions, one needs to
distinguish between different types of actions of subgroups of \glgroup{N}{\complex}. Note that any $G\subset\glgroup{N}{\complex}$ comes equipped with an action on
$\complex^{N}$. Hence the properties of the action of $G$ can be said to be
properties of $G$ (or $\bar{G}$) itself.

\begin{defin}
The group $G$ is \emph{irreducible} if for any non-zero vector $x\in\complex^{N}$,
the $G$-orbit of $x$ spans $\complex^{N}$.
\end{defin}
\begin{defin}
The group $G$ is \emph{primitive} if there does not exist a decomposition of
$\complex^{N}$ into a direct sum of proper linear subspaces, such that $G$
permutes the subspaces. If such a decomposition does exist, then the action is
\emph{imprimitive}.
\end{defin}
\begin{defin}
An imprimitive group $G$ is \emph{monomial}, if the action of $G$
is induced from a $1$-dimensional action of some subgroup $G'\subseteq G$. In
other words, $G=D\rtimes T$, where $D$ consists of diagonal matrices, and
$T\subseteq\symgroup{N}$ acts by permuting the basis. In order to avoid
repetitions, throughout this paper ``monomial'' will be used to mean ``irreducible
monomial'' (in particular, $T$ is assumed to be transitive).
\end{defin}

From the more geometric point of view, it is more interesting to look at the
variety $\complex^{n}/G$. Note that since this is a quotient singularity, it
is a Kawamata log terminal singularity.

\begin{defin}
\label{definition:exceptional} Let $(V\ni O)$ be a~germ of a Kawamata log
terminal singularity. The singularity is~said to be \emph{exceptional} if for every
effective \rationals-divisor $D_{V}$ on the~variety $V$, such that the~log pair $(V,D_{V})$ is log canonical,
there exists at most one exceptional divisor over the point $O$ with discrepancy $-1$ with respect to the pair $(V,D_{V})$.
\end{defin}

\begin{prop}[see~{\cite[Proposition~2.1]{Prokhorov00}}]\label{prop:Exc-prim}
Let $G\subset\slgroup{N}{\complex}$ be a finite subgroup that induces an exceptional
singularity. Then $G$ is primitive.
\end{prop}
\begin{cor}
For any given $N$, only finitely many finite subgroups of \slgroup{N}{\complex}
induce exceptional singularities.
\end{cor}
\pf Immediate by Proposition~\ref{prop:Exc-prim} and Jordan's theorem
(see, for example, \cite{Frobenius11}).\qed

One can define a larger class of singularities by looking at a special class of birational morphisms:
\begin{thm}[see~{\cite[Theorem~3.7]{Cheltsov-Shramov09}}]\label{thm:plt}
Let $(V\ni O)$ be a~germ of a~Kawamata log terminal singularity. Then there
exists a birational morphism \deffunname{\pi}{W}{V} such that the following
hypotheses are satisfied:
\begin{itemize}
\item the~exceptional locus of $\pi$ consists of one irreducible divisor
$E$ such that $O\in\pi(E)$,%
\item the~log pair $(W,E)$ has purely log terminal singularities.
\item the~divisor $-E$ is a~$\pi$-ample \rationals-Cartier divisor.
\end{itemize}
\end{thm}
\pf The result follows from~\cite{BCHM},
\cite[Theorem~1.5]{kudryavtsev2001}
and~\cite[Proposition~2.9]{Prokhorov00blowups} \qed
\begin{defin}
Let $(V\ni O)$ be a~germ of a~Kawamata log terminal singularity,
and \deffunname{\pi}{W}{V} be~a~birational morphism satisfying the conditions
of Theorem~\ref{thm:plt}. Then $\pi$ is a~\emph{plt blow-up} of
the~singularity.
\end{defin}
\begin{defin}
We say that the singularity $(V\ni O)$ is \emph{weakly-exceptional} if it has a unique plt blow-up.
\end{defin}
\begin{lemma}[see~{\cite[Theorem~4.9]{Prokhorov00blowups}}]
If $(V\ni O)$ is exceptional, then $(V\ni O)$ is weakly-exceptional.
\end{lemma}
\begin{thm}\label{thm:WExc-trans}
Let $G\subset\slgroup{N}{\complex}$ be a finite subgroup that induces a
weakly-exceptional singularity. Then $G$ is irreducible.
\end{thm}
\pf The argument is similar to that in~\cite[Proposition~2.1]{Prokhorov00}. \qed

This result will be used frequently throughout this paper, since for any finite
group $G$, it allows to severely limit the number of actions worth checking.

\begin{eg}
One can consider the case $N=2$. It is a well-known fact that the finite
subgroups of \Aut{\proj{1}} are (up to conjugation): cyclic (\cyclgroup{n},
$n\geq1$), dihedral (\dihgroup{2n}, $n\geq2$) or polyhedral (\altgroup{4},
\symgroup{4}, \altgroup{5}), with lifts to \slgroup{2}{\complex} being their
central extensions by \cyclgroup{2} (except the case of odd cyclic groups, where the central extension is not necessary). Cyclic groups induce reducible actions,
polyhedral groups induce primitive actions, and the dihedral groups induce
monomial actions.

On the geometric side, rephrasing~\cite[Section~5.2.3]{Shokurov92-Eng} implies
that \altgroup{4}, \symgroup{4} and \altgroup{5} give rise to
exceptional singularities of types $E_{6}$, $E_{7}$ and $E_{8}$ respectively,
the dihedral groups \dihgroup{n} give rise to weakly-exceptional, but not
exceptional singularities of type $D_{n+2}$ and the cyclic groups \cyclgroup{n}
give rise to singularities of type $A_{n-1}$, that are not weakly-exceptional.
The concept of a weakly-exceptional singularity was originally formed as a
generalisation of the singularities of types $D_{n}$ and $E_{n}$.
\end{eg}

Since this example fully settles the $N=2$ case, assume for the rest of this paper
that $N\geq3$. In the dimensions $N\leq6$, the exceptional singularities have
been fully classified in \cite{Cheltsov-Shramov09}, \cite{Cheltsov-Shramov10}
and \cite{Markushevich-Prokhorov}. The aim of this paper is to obtain a list of
group actions that induce weakly-exceptional singularities in dimension $3$
and $4$. This will be done by applying the exceptionality criteria obtained by
Y.~Prokhorov, I.~Cheltsov and C.~Shramov in~\cite{Cheltsov-Shramov09}
and~\cite{Prokhorov00}, using the language of Tian's alpha-invariant introduced
in~\cite{Tian87} and~\cite{Tian-Yau87} by G.~Tian and S-T.~Yau.

For dimension $N=3$, the classification of finite subgroups of
\slgroup{3}{\complex} is a classical result by G.~Miller, H.~Blichfeldt and
L.~Dickson (see~\cite{MDB}). A modern exposition of this result has been made
by S.~Yau and Y.~Yu (see~\cite[Theorem~A]{Yau-Yu}). 
A result by D.~Markushevich and Y.~Prokhorov says:
\begin{thm}[see~\cite{Markushevich-Prokhorov}]
The group $G$ induces an exceptional singularity if and only if $\bar{G}$ is
isomorphic to \altgroup{6}, Klein's simple group $\mathbb{K}_{168}$ of size
$168$,
Hessian group $\mathbb{H}_{648}$ of size $648$ or its normal subgroup $F_{216}$ of size $216$.
\end{thm}
The first main result of this paper extends this by:
\begin{thm}\label{sl3:thm:main}
Let $G\subset\slgroup{3}{\complex}$ be a finite subgroup. Then $G$ induces
a weakly-exceptional but not an exceptional singularity if and only if one of the following holds:
\begin{itemize}
 \item $G$ is a monomial group, and $\bar{G}$ is not isomorphic to
$\left(\cyclgroup{2}\right)^{2}\rtimes\cyclgroup{3}$ or
$\left(\cyclgroup{2}\right)^{2}\rtimes\symgroup{3}$.
 \item $G$ is isomorphic to the normal subgroup $E_{108}\lhd F_{216}$ of size
$108$.
\end{itemize}
\end{thm}
Section~\ref{sect:3D} of this paper is devoted to proving this theorem.

Since $4$ is not a prime number, \slgroup{4}{\complex} contains significantly
more finite subgroups than \slgroup{3}{\complex} does. The list of finite
primitive subgroups of \slgroup{4}{\complex} is a classical result, that can be found in H.~Blichfeldt's
book (see~\cite[Chapter~VII]{Blichfeldt}). The list of irreducible imprimitive finite
subgroups can be obtained from papers by D.~Flannery (see~\cite{Flannery}) and
B.~H\"{o}fling (see~\cite{Hofling}).

The second main result of this paper is designed to supplement the list of finite subgroups of
\slgroup{4}{\complex} inducing exceptional singularities, that can be found in a paper by I.~Cheltsov and C.~Shramov (see~\cite{Cheltsov-Shramov09}).
In order to avoid reproducing the tables of groups from the works mentioned above, the paper will instead produce a list of
irreducible group actions that give rise to singularities that are \emph{not} weakly-exceptional. Therefore, the second main result of this paper is:
\begin{thm}\label{sl4:thm:main}
If $G\subset\slgroup{4}{\complex}$ is a finite subgroup whose action
induces a singularity that is
not weakly-exceptional, then either $G$ is not irreducible, or $\bar{G}$ must be
conjugate to one of (using the notation described in
Section~\ref{sect:quadric-action} where appropriate):
\begin{itemize}
\item A primitive group that is one of:
\begin{itemize}
\item \symgroup{5}, \altgroup{5}.

\item 9 groups $\left(H_{1},H_{1},H_{2},H_{2}\right)\cong H_{1}\times
H_{2}$, for different choices of
$H_{1},H_{2}\in\defsetshort{\altgroup{4},\symgroup{4},\altgroup{5}}$.

 \item $\left(\altgroup{4}\times\altgroup{4}\right)\rtimes\cyclgroup{2}$,
$\left(\symgroup{4}\times\symgroup{4}\right)\rtimes\cyclgroup{2}$,
$\left(\altgroup{5}\times\altgroup{5}\right)\rtimes\cyclgroup{2}$.

 \item
$\frac{1}{2}\left[\symgroup{4}\times\symgroup{4}\right]
\cong\left(\altgroup{4}\times\altgroup{4}\right)\rtimes\cyclgroup{2}
\cong\left(\symgroup{4},\altgroup{4},\symgroup{4},\altgroup{4}\right)$.

 \item
$\left(\left(\altgroup{4}\times\altgroup{4}
\right)\rtimes\cyclgroup{2}\right)\rtimes\cyclgroup{2}
\cong\left(\symgroup{4},\altgroup{4},\symgroup{4},\altgroup{4}
\right)\rtimes\cyclgroup{2}$.

\end{itemize}
\item An imprimitive non-monomial group that is one of:
\begin{itemize}
\item 3 families of groups $\dihgroup{2m}\times H_{2}$, where
$H_{2}\in\defsetshort{\altgroup{4},\symgroup{4},\altgroup{5}}$.

 \item
$\frac{1}{2}\left[\dihgroup{2m}\times\symgroup{4}\right]
\cong\left(\cyclgroup{m}\times\altgroup{4}\right)\rtimes\cyclgroup{2}
\cong\left(\dihgroup{2m},\cyclgroup{m},\symgroup{4},\altgroup{4}\right)$.

 \item
$\frac{1}{2}\left[\dihgroup{4m}\times\symgroup{4}\right]
\cong\left(\dihgroup{2m}\times\altgroup{4}\right)\rtimes\cyclgroup{2}
\cong\left(\dihgroup{4m},\dihgroup{2m},\symgroup{4},\altgroup{4}\right)$
for $m\geq2$.

 \item
$\frac{1}{6}\left[\dihgroup{6m}\times\symgroup{4}\right]
\cong\left(\cyclgroup{m}\times \mathbb{V}_{4}\right)\rtimes\symgroup{3}
\cong\left(\dihgroup{6m},\cyclgroup{m},\symgroup{4},\mathbb{V}_{4}\right)$.

\end{itemize}
\item A monomial group that is one of:
\begin{itemize}
 \item
$\left(\left(\cyclgroup{3}\right)^{3}\rtimes\cyclgroup{2}
\right)\rtimes\cyclgroup{2}$, acting as shown in
Lemma~\ref{lemma:sl4:cubic:smooth}.
 \item $\left(\cyclgroup{3}\right)^{3}\rtimes\symgroup{4}$, acting as shown in
Lemma~\ref{lemma:sl4:cubic:smooth}.

 \item \altgroup{4} or \symgroup{4}, acting as twisted diagonal groups $\left(H_{1},1,H_{1},1\right)_{\alpha}$.

 \item $\dihgroup{2m}\times\dihgroup{2n}$ ($m,n\geq2$).

 \item
$\left(\dihgroup{2n}\times\dihgroup{2n}\right)\rtimes\cyclgroup{2}$ ($m,n\geq2$).

 \item $\dihgroup{4n}\rtimes\cyclgroup{2}$ ($n\geq2$).
 \item $\altgroup{4}\rtimes\cyclgroup{2}$, $\symgroup{4}\rtimes\cyclgroup{2}$.

 \item
$\frac{1}{6}\left[\symgroup{4}\times\symgroup{4}\right]
\cong \left(\mathbb{V}_{4}\times \mathbb{V}_{4}\right)\rtimes\symgroup{3}
\cong\left(\symgroup{4},\mathbb{V}_{4},\symgroup{4},\mathbb{V}_{4}\right)$.

 \item
$\frac{1}{3}\left[\altgroup{4}\times\altgroup{4}\right]
\cong\left(\mathbb{V}_{4}\times \mathbb{V}_{4}\right)\rtimes\cyclgroup{3}
\cong\left(\altgroup{4},\mathbb{V}_{4},\altgroup{4},\mathbb{V}_{4}\right)$.

 \item
$\frac{1}{2}\left[\dihgroup{2m}\times\dihgroup{4n}\right]
\cong\left(\cyclgroup{m}\times\dihgroup{2n}\right)\rtimes\cyclgroup{2}
\cong\left(\dihgroup{2m},\cyclgroup{m},\dihgroup{4n},\dihgroup{2n}\right)$
($m,n\geq2$).

 \item
$\frac{1}{4}\left[\dihgroup{4m}\times\dihgroup{4n}\right]_{\alpha}
\cong\left(\cyclgroup{m}\times \cyclgroup{n}\right)\rtimes\dihgroup{4}
\cong\left(\dihgroup{4m},\cyclgroup{m},\dihgroup{4n},\cyclgroup{n}\right)_{
\alpha }$ (where
$\alpha(b)=a_{2n}^{n}$, $\alpha(a_{2m}^{m})=b$).

 \item
$\frac{1}{2}\left[\dihgroup{4m}\times\dihgroup{4n}\right]
\cong\left(\dihgroup{2m}\times\dihgroup{2n}\right)\rtimes\cyclgroup{2}
\cong\left(\dihgroup{4m},\dihgroup{2m},\dihgroup{4n},\dihgroup{2n}\right)$
($m,n\geq2$).

 \item
$\left(\left(\mathbb{V}_{4}\times
\mathbb{V}_{4}\right)\rtimes\symgroup{3}\right)\rtimes\cyclgroup{2}
\cong\left(\symgroup{4},\mathbb{V}_{4},\symgroup{4},\mathbb{V}_{4}
\right)\rtimes\cyclgroup{2}$.

 \item
$\left(\left(\mathbb{V}_{4}\times
\mathbb{V}_{4}\right)\rtimes\cyclgroup{3}\right)\rtimes\cyclgroup{2}
\cong\left(\altgroup{4},\mathbb{V}_{4},\altgroup{4},\mathbb{V}_{4}
\right)\rtimes\cyclgroup{2}$.

 \item
$\left(\left(\cyclgroup{m}\times\cyclgroup{m}\right)\rtimes\dihgroup{4}
\right)\rtimes\cyclgroup{2}
\cong\left(\dihgroup{4m},\cyclgroup{m},\dihgroup{4m},
\cyclgroup{m}\right)_{\alpha}\rtimes\cyclgroup{2}$
(where
$\alpha(b)=a_{2m}^{m}$, $\alpha(a_{2m}^{m})=b$).

 \item
$\left(\left(\dihgroup{2m}\times\dihgroup{2m}\right)\rtimes\cyclgroup{2}
\right)\rtimes\cyclgroup{2}
\cong\left(\dihgroup{4m},\dihgroup{2m},\dihgroup{4m},\dihgroup{2m}
\right)\rtimes\cyclgroup{2}$ for $m\geq2$.

\end{itemize}
\end{itemize}
\end{thm}
Section~\ref{sect:4D} of this paper is devoted to proving this result.
\begin{rmk}
Note that the list above gives explicit conjugacy classes of most of the
subgroups. In most cases, the irreducible subgroup is uniquely defined (up
to conjugacy) by its isomorphism class.
\end{rmk}

The results in this paper can be applied in birational geometry,
including the problem of conjugacy of subgroups of higher-dimensional Cremona
groups. For more details, see~\cite{Cheltsov08delPezzo}.

\section{Preliminaries}
\subsection{Notation}\label{notation} In this section, some standard notation used throughout this paper is defined:
\begin{itemize}
 \item \cyclgroup{n} is the cyclic group of size $n$.
 \item \dihgroup{2n} is the dihedral group of size $2n$.
 \item \symgroup{n} is the permutation group of a set of $n$ elements.
 \item \altgroup{n} is the alternating group, an index 2 subgroup of \symgroup{n}.
 \item $\bcyclgroup{n}, \bdihgroup{2n}, \baltgroup{n}, \bsymgroup{n}\subset\slgroup{2}{\complex}$ are the
binary versions of the relevant groups, i.e.\ their central extensions by
\cyclgroup{2} (see, for example, \cite[Section~4.4]{Springer77}).
Since the generators of these groups will be referred to heavily at the end
of Section~\ref{sect:4D}, fix the presentations of these groups as:
 \begin{itemize}
   \item $\bcyclgroup{n}=\defsetspanlong{a_{n}}{a_{n}^{2n}=1}$
   \item $\bdihgroup{2n}=\defsetspanlong{a_{n},b}{a_{n}^{n}=b^{2}, b^{4}=1, ba_{n}b^{-1}=a_{n}^{-1}}$
   \item $\baltgroup{4}=\defsetspan{[12][34], [14][23], [123]}$, where the basis is chosen so that $[12][34]$ preserves the two lines spanned by the basis vectors, and $[14][23]$ swaps them. Extend this presentation to that of $\bsymgroup{4}=\defsetspan{[12][34], [14][23], [123], [34]}$.
   \item $\baltgroup{5}=\defsetspan{[12345], [12][34]}$, where the basis is chosen so that $[12345]$ is diagonal.
 \end{itemize}
The last $3$ groups are central extensions of permutation groups, and their generators are intentionally named to identify the permutations they correspond to. The relations come from the relations between the permutations.
 \item $\mathbb{V}_{4}$ is the Klein group of size $4$.
 \item $\zeta_{n}$ is a primitive $n$-th root of unity. Whenever two such are
used in defining generators of the same group, the choice is assumed to be
consistent, i.e.\ $\zeta_{ab}^{a}=\zeta_{b}$ for $a,b\in\integers$.
 \item Write ``extension of $G$ by scalar elements'' (for a finite group $G$) to mean $H\rtimes G\subset\slgroup{N}{\complex}$, where $H$ is a subgroup of the centre of \slgroup{N}{\complex} (i.e.\ consists of scalar matrices). Note that unless $N$ is prime, there are several non-trivial possibilities for $H$ that give $H\rtimes G$ the same image under the natural projection to \pglgroup{N}{\complex}.
\end{itemize}
Note that the groups $\mathbb{V}_{4}$, \dihgroup{4} and
$\cyclgroup{2}\times\cyclgroup{2}$
are isomorphic. The notation used for this group will denote the context the group is considered in: $\mathbb{V}_{4}$ will be used whenever it is considered as a subgroup of \altgroup{4} or \symgroup{4}, and one of the
other two will be used whenever it is considered on its own.

\subsection{Group actions on a smooth quadric surface in \proj{3}}\label{sect:quadric-action}
In this section some notation for a specific type of group action on \proj{3} will be presented. This is a general form for an action that preserves a smooth quadric surface in \proj{3} (which can be taken to be $\proj{1}\times\proj{1}$ embedded via the Segre embedding). Some of the notation was taken from this action's description in~\cite[Section~4.3]{D-I}, with some additions tailored to describing individual members of families of related groups.

It is possible to present \proj{3} as the set of non-zero $2\times2$ matrices modulo the scalar ones.
From now on, consider the ``matrix form'' of \proj{3} to be:
\[\left(x:y:u:v\right)\leadsto\left(\begin{array}{cc}x&y\\u&v\end{array}\right)\]

Let $S$ be the image of $\proj{1}\times\proj{1}$ under the Segre embedding into $\proj{3}$. Then one can assume $S$ is the zero set of the determinant of the matrix form of \proj{3}.

Let $\bar{G}$ be a finite group acting faithfully on $\proj{1}\times\proj{1}$.
This variety has exactly two rulings, which $\bar{G}$ can either preserve or
interchange. Consider the exact sequence
\[1\longrightarrow H\longrightarrow \bar{G}\longrightarrow
\symgroup{2}\]
where \symgroup{2} permutes the two rulings. Then $H\unlhd\bar{G}$ is the
maximal subgroup that preserves the ruling, and either $\bar{G}=H$ or $\bar{G}$ is generated by $H$ and an element $\sigma\tau$, where $\tau$ is the involution interchanging the two rulings, and $\sigma$ is some automorphism of $S$ preserving the ruling, with $\sigma^{2}\in H$ (see~\cite[Theorem~4.9]{D-I}).

Let \deffunname{\pi_{i}}{H}{H_{i}} be the projections of $H$ on the
two components of $\proj{1}\times\proj{1}$. Then have two more short exact
sequences
\[1\longrightarrow K_{2}\longrightarrow H\longrightarrow H_{1}\longrightarrow1\]
\[1\longrightarrow K_{1}\longrightarrow H\longrightarrow H_{2}\longrightarrow1\]
It is clear that $K_{1}\cap K_{2}=\defsetshort{1}$. Therefore, for $i\neq j$
($i,j\in\defsetshort{1,2}$),
\[
K_{i}\cong\hat{K}_{i}\defby\defset{kK_{j}}{k\in K_{i}}\unlhd H/K_{j}=H_{i}
\]
In this notation, $H_{1}/\hat{K}_{1}\cong H/\left(K_{1}K_{2}\right)\cong
H_{2}/\hat{K}_{2}$, so the group can be defined completely by $\left(H_{1},K_{1},H_{2},K_{2}\right)_{\alpha}$, where $\alpha$ is an isomorphism
\deffun{H_{1}/\hat{K}_{1}}{H_{2}/\hat{K}_{2}}. In return, if $H$ is known, one
can reconstruct $\alpha$ by making it map
\defmap{\pi_{1}(h)\hat{K}_{1}}{\pi_{2}(h)\hat{K}_{2}} ($\forall h\in H$).

In the matrix form described above, $H_{1}$ acts on \proj{3} by left matrix
multiplication, and $H_{2}$ acts by transposed right matrix multiplication.
The involution switching the two rulings of $\proj{1}\times\proj{1}$ corresponds
to transposing the matrix. If $h_{1}\in H_{1}$ acts on the first component of
$\proj{1}\times\proj{1}$ and $h_{2}\in H_{2}$ acts on the second one, then
write $\quadricaction{h_{1}}{h_{2}}$ to denote this action. Explicitly, for any
$2\times2$ matrices $A$ and $B$,
\[
\quadricaction{A}{B}\left(
\left(\begin{array}{cc}x&y\\u&v\end{array}\right)\right) =
A\left(\begin{array}{cc}x&y\\u&v\end{array}\right)B^{T}
\]
Furthermore, interchanging the rulings of $\proj{1}\times\proj{1}$ corresponds
to transposing the matrix form of \proj{3}. It is worth noting that if $\bar{G}\neq H$,
then conjugation by $\sigma\tau$ provides isomorphisms $H_{1}\cong H_{2}$ and
$K_{1}\cong K_{2}$.

The following notation may sound somewhat unnecessary, but it will assist in avoiding describing explicit generators of groups later on.
For some groups, one can say some of its elements are of some special
``type'', e.g.\ the order $2$ element in \dihgroup{2n} (n odd) or the
$3$-cycles in \altgroup{4}. Given a finite subgroup
$\left(H_{1},K_{1},H_{2},K_{2}\right)_{\alpha}$ with
$H_{1},H_{2}\subset\slgroup{2}{\complex}$, say elements of $H_{1}$ and $H_{2}$
of some fixed type are ``coupled'', if $\forall\quadricaction{h}{h'}\in H$ with
$h$ an element of this type, then $h'$ must be either an element of the same
type or a product of such an element and an element of a different type.
Otherwise say that elements of this type are ``not coupled''. For example, if
$H_{1}\cong H_{2}\cong\dihgroup{10}$, then order $2$ elements $b$ are not
coupled in $G$ if $G$ contains an element $\quadricaction{a^{k}}{ba^{l}}$ for
some $k,l\in\integers$, where $a$ is an order $5$ element. Otherwise they are
coupled.

\subsection{Exceptionality criteria}\label{sect:criteria}
In order to determine the exceptionality of a given singularity, it is more
useful to have a more computable criterion than that given by the definitions
of exceptionality and weak exceptionality. To discuss this criterion, some definitions
need to be made.
\begin{defin}
Let $X$ be a smooth Fano variety (see~\cite{Iskovskikh-Prokhorov99}) of
dimension $n$, and let $g=\left(g_{ij}\right)$ be a K\"{a}hler metric, such that
\[\omega=\frac{\sqrt{-1}}{2\pi}\sum{g_{ij}\mbox{d}z_{i}\wedge\mbox{d}\bar{z}_{j}
}\in c_{1}\left(X\right) \]
Let $\bar{G}\subseteq\Aut{X}$ be a compact subgroup, such that $g$ is
$\bar{G}$-invariant.
Let $\mbox{P}_{\bar{G}}\left(X,g\right)$ be the set of $C^{2}$ smooth
$\bar{G}$-invariant functions, such that $\forall
\phi\in\mbox{P}_{\bar{G}}\left(X,g\right)$,
\[\omega+\frac{\sqrt{-1}}{2\pi}\partial\bar{\partial}\phi>0\]
and $\sup_{X}{\phi}=0$. Then the \emph{$\bar{G}$-invariant $\alpha$-invariant}
of
$X$ is
\[
\alpha_{\bar{G}}\left(X\right)=\sup{\defset{\lambda\in\rationals
}{\begin{array}{c}\exists C\in\reals,\ \mbox{such that}\
\forall\phi\in\mbox{P}_{\bar{G}}\left(X,g\right),\\
\int_{X}{e^{-\lambda\phi}\omega^{n}}\leq C \end{array}}}
\]
where $n$ is the dimension of $X$.
\end{defin}
The number $\alpha_{\bar{G}}\left(X\right)$ was introduced in~\cite{Tian87}
and~\cite{Tian-Yau87}.

\begin{defin}
Let $X$ be a variety with at most Kawamata log terminal singularities
(see~\cite[Definition~3.5]{Kollar97}) and $D$ an effective $\rationals$-divisor
on $X$. Let $Z\subseteq X$ be a closed non-empty subvariety. Then the \emph{log
canonical threshold} of $D$ along $Z$ is
\[c_{Z}\left(X,D\right)=\sup{\defset{\lambda\in\rationals}{\mbox{the pair}\ 
\left(X,\lambda D\right)\ \mbox{is log canonical along}\ Z}}\]
To simplify notation, write $c_{X}\left(X,D\right)=c\left(X,D\right)$.
\end{defin}
There exists an equivalent complex
analytic definition of the log canonical threshold:
\begin{prop}[see~{\cite[Proposition~8.2]{Kollar97}}]
Let $X$ be a smooth complex variety, $Z$ a closed non-empty subscheme of $X$,
and $f$ a non-zero regular function on $X$. Then
\[
c_{Z}\left(X,\defsetshort{f=0}\right)=\sup{\defset{\lambda\in\rationals}{
\abs{f}^{-\lambda}\ 
\mbox{is locally}\ L^{2}\ \mbox{near}\ Z }}
\]
\end{prop}

\begin{defin}
Let $G$ be a~finite subgroup of \glgroup{N}{\complex},
where $N\geq2$, and let $\bar{G}$ be its image under the natural projection
into \pglgroup{N}{\complex}. Then the \emph{global $\bar{G}$-invariant log
canonical
threshold} of $\proj{N-1}$ is:
\[
\mathrm{lct}\left(\proj{N-1},\bar{G}\right)=\inf{\defset{c\left(\proj{N-1},
D\right)}{\begin{array}{l}
D\ \mbox{is a $\bar{G}$-invariant effective}\\
\mbox{\rationals-divisor on \proj{N-1}, such that}\\
D\sim_{\rationals}-\canclass{\proj{N-1}}
\end{array}}} 
\]
\end{defin}
\begin{rmk}[see~{\cite[Theorem~A.3]{Cheltsov-Shramov08-appendix}}]
$\mathrm{lct}\left(\proj{N-1},\bar{G}\right)=\alpha_{\bar{G}}\left(\proj{N-1}
\right)$
\end{rmk}

\begin{thm}[see {\cite[Theorem~3.15]{Cheltsov-Shramov09}}]
The~singularity $\complex^{N}/G$ is weakly-exceptional $\iff$
$\mathrm{lct}(\proj{N-1},\bar{G})\geq1$.
\end{thm}

A similar condition is often necessary in order to compute conjugacy classes in
higher-dimensional Cremona groups. For details,
see~\cite{Cheltsov08delPezzo}.

\begin{thm}[see {\cite[Theorem~3.18]{Cheltsov-Shramov09}}]
\label{thm:criterion:3D} Suppose that
$G\subset\slgroup{3}{\complex}$ is a finite group and $\bar{G}$ is its image under
the natural projection into \pglgroup{3}{\complex}. Then the~following are
equivalent:
\begin{itemize}
\item the~inequality $\mathrm{lct}(\proj{2},\bar{G})\geq1$ holds,
\item the~group $G$ does not have semi-invariants of degree at most $2$.
\end{itemize}
\end{thm}

\begin{thm}[see {\cite[Theorem~4.1]{Cheltsov-Shramov09}}]
Suppose that $G\subset\slgroup{4}{\complex}$ is a finite group and $\bar{G}$ is its
image under the natural projection into \pglgroup{4}{\complex}. Then
the~inequality $\mathrm{lct}(\proj{3},\bar{G})\geq1$ holds if and only if
the~following conditions are satisfied:
\begin{itemize}
\item the~group $G$ is irreducible,
\item the~group $G$ does not have semi-invariants of degree at most $3$,
\item there is no $\bar{G}$-invariant smooth rational cubic curve in
$\proj{3}$.
\end{itemize}
\end{thm}
Equivalently, the same criterion can be stated as:
\begin{thm}\label{thm:criterion:4D}
Suppose that $G\subset\slgroup{4}{\complex}$ is a finite group and $\bar{G}$ is its
image under the natural projection into \pglgroup{4}{\complex}. Then
the~inequality $\mathrm{lct}(\proj{3},\bar{G})\geq1$ holds if and only if
the~following conditions are satisfied:
\begin{itemize}
\item the~group $G$ is irreducible,
\item the~group $G$ does not have semi-invariants of degree at most $3$,
\item $G$ is not a central extensions of \altgroup{5} acting on
$\complex^{4}$ as the third symmetric power of its irreducible $2$-dimensional
representation.
\end{itemize}
\end{thm}
\pf
Assume $\bar{G}\subset\pglgroup{4}{\complex}$ is a finite irreducible subgroup that
does not fix any quadric or cubic surface, and
$C\subset\proj{3}$ is a smooth rational $\bar{G}$-invariant cubic curve.

One can assume $C$ is the image of
\[\defmapname{\mathcal{C}}{\left(x:y\right)\in\proj{1}
}{\left(x^{3}:x^{2}y:xy^{2}:y^{3}\right)\in\proj{3}}\]
This easily implies (by Proposition~\ref{prop:SES}) that $\bar{G}$ must be
isomorphic to one of the finite automorphism
groups of \proj{1}, with its action induced by the action of \proj{1} via
$\mathcal{C}$. This means $\bar{G}$ must be one of the following:
\begin{itemize}
 \item Cyclic or dihedral group: neither these groups nor their central
extensions have an irreducible
$4$-dimensional representation. Therefore such a $G$ would not be irreducible.
 \item \altgroup{4} or \symgroup{4}: The induced actions of these two groups
preserve \defset{\left(x:y:u:v\right)\in\proj{3}}{xy-uv=0}, which is a smooth
quadric surface.
 \item \altgroup{5}: This action is primitive. From the embedding of $C$ into
\proj{3}, it is easy to see that the action on $\complex^{4}$ must correspond to either
the third symmetric power of the irreducible $2$-dimensional representation of
\baltgroup{5} or its central extension.
\end{itemize}
\qed
\begin{rmk}
The last condition of Theorem~\ref{thm:criterion:4D} is necessary, since this
action is irreducible and does not
preserve any
projective surfaces of degree $2$ or $3$ (can be checked directly or seen in, for example,
\cite[Proof of Lemma~4.9]{Cheltsov-Shramov09}).
\end{rmk}

\subsection{Properties of possible invariants}\label{prelim:invars}
Let $\bar{G}\subset\pglgroup{N}{\complex}$ be a finite subgroup, and $G$ be its
lift
to \slgroup{N}{\complex}. Assume that $S\subset\proj{N-1}$ is a proper
$\bar{G}$-invariant subvariety. Then there exists a natural homomorphism
\deffunname{\pi_{S}}{G}{\Aut{S}} giving a short exact sequence
\[0\longrightarrow G_{0}\longrightarrow G\stackrel{\pi_{S}}{\longrightarrow}
G_{S}\longrightarrow0\]
where $G_{S}=\pi_{S}\left(G\right)$, $G_{0}=\ker{\pi_{S}}$.

\begin{prop}\label{prop:SES}
Let $G_{0}\subset\slgroup{N}{\complex}$ be a finite subgroup, and
$\bar{G}_{0}$ the image of its natural embedding into \pglgroup{N}{\complex}.
Let $S\subset\proj{N-1}$ be a subvariety that is not contained in
the union of any two proper linear subspaces of \proj{N-1},
and $\bar{G}_{0}$ fixes $S$ point-wise. Then $\bar{G}_{0}$ is trivial.
\end{prop}

\pf  Pick $g\in G_{0}$. Then $\defsetspan{g}$ is a finite
abelian group, and so (in some basis for $\complex^{N}$) consists of diagonal
matrices. Let $\bar{g}$ be the image of $g$ under the natural projection
into \pglgroup{N}{\complex}.

Let $\bar{e}_{1},\ldots,\bar{e}_{N}\in S$ be distinct points, such that 
their lifts $e_{1},\ldots,e_{N}$ to $\complex^{N}$ span all of $\complex^{N}$
(these exist, since $S$ is not contained in a proper linear subspace of
\proj{N-1}). Then $e_{i}$ are eigenvectors of $g$, and let $\lambda_{i}$ be the
corresponding eigenvalues. Reordering $e_{i}$-s if necessary, let
$1\leq m\in\integers$ be such that $\lambda_{1}=\ldots=\lambda_{m}$ and
$\lambda_{n}\neq\lambda_{m}$ $\forall m<n\leq N$.

Assume $m<N$. Then take $A\subset\complex^{N}$ to be the linear subspace
spanned by $e_{1},\ldots,e_{m}$ and $B\subset\complex^{N}$ to be the
linear subspace spanned by $e_{m+1},\ldots,e_{N}$. Let $\bar{A},\bar{B}$ be
their natural projections into \proj{N-1}. These are proper linear subspaces, so
$\exists \bar{p}\in S\setminus\left(\bar{A}\cup\bar{B}\right)$.

This means there are at least $1+N-m$ distinct linear eigenspaces for $g$ not
contained in $A$, at least one of which is not contained in $B$ either.
Therefore must have $\lambda_{n}=\lambda_{m}$ for some $m<n\leq N$,
contradicting the choice of $m$.

This means $m=N$, and so $g$ is a scalar matrix.\qed

\begin{cor}
Let $G\subset\slgroup{N}{\complex}$ and let $\bar{G}$ be its natural projection into
\pglgroup{N}{\complex}. Let $S\subset\proj{N-1}$ be a $\bar{G}$-invariant
subvariety, that is not contained in the union of any two proper linear
subspaces of \proj{N-1}. Let \deffunname{\pi_{s}}{\bar{G}}{\Aut{S}} be the
natural homomorphism. Then $\ker{\pi_{S}}=0$.
\end{cor}

\begin{rmk}
If $S\subset\proj{N-1}$ is an irreducible surface, then either it is contained
in an $(N-2)$-dimensional linear subspace of \proj{N-1} or it is not contained
in the union of any two proper linear subspaces of \proj{N-1}.
\end{rmk}

Now assume further that $G$ is irreducible. This implies the following
restrictions:
\begin{rmk}\label{rmk:excluding-1} In the notation above, the following hold:
\begin{itemize}
 \item $\bar{G}$ cannot be cyclic, as in that case $G$ is abelian (as it is a
central extension by \emph{scalar} elements), and all irreducible
representations of abelian groups are $1$-dimensional.
 \item $S$ cannot be contained in a proper linear subspace of \proj{N-1}.
\end{itemize}
\end{rmk}

\begin{lemma}\label{lemma:singNo}
In the notation above, let $S$ have an isolated singularity of some given type
(e.g.\ $A_{n}$, $D_{n}$, etc.). Then $S$ has at least $N$ isolated
singularities of this type.
\end{lemma}
\pf Let $S$ have exactly $k$ ($1\leq k<N$) singularities of this type. Then they form a $\bar{G}$-invariant set of $k$ points, giving a $k$-dimensional $\bar{G}$-invariant subspace of \proj{N}. This is impossible, as $G$ is irreducible.\qed

\section{Three-dimensional case}\label{sect:3D}
Assume $N=3$ throughout this section. This section is devoted to proving
Theorem~\ref{sl3:thm:main}.

The list of all finite subgroups of \slgroup{3}{\complex} is given by:
\begin{prop}[see~{\cite[Theorem~A]{Yau-Yu}}]\label{sl3:prop:classification}
Define the following matrices:
\[\begin{array}{lll}
S=\left(\begin{array}{ccc}1&0&0\\0&\omega&0\\0&0&\omega^{2}\end{array}\right)&
T=\left(\begin{array}{ccc}0&1&0\\0&0&1\\1&0&0\end{array}\right)&
W=\left(\begin{array}{ccc}\omega&0&0\\0&\omega&0\\0&0&\omega\end{array}
\right)\\
&&\\
U=\left(\begin{array}{ccc}\epsilon&0&0\\0&\epsilon&0\\0&0&\epsilon\omega
\end{array}\right)&
Q=\left(\begin{array}{ccc}a&0&0\\0&0&b\\0&c&0\end{array}\right)&
V=\frac{1}{\sqrt{-3}}\left(\begin{array}{ccc}1&1&1\\1&\omega&\omega^{2}
\\1&\omega^{2}&\omega\end{array}\right)\\
\end{array}\]
where $\omega=e^{2\pi i/3}$, $\epsilon^{3}=\omega^{2}$ and
$a,b,c\in\complex$ are chosen arbitrarily, as long as $abc=-1$ and
$Q$ generates a finite group.

Up to conjugacy, any finite subgroup of \slgroup{3}{\complex} belongs to one of
the following types:
\begin{enumerate}[(A)]
 \item\label{gptype:A} Diagonal abelian group.
 \item\label{gptype:B} Group isomorphic to an irreducible finite subgroups of
\glgroup{2}{\complex} and not conjugate to a group of type~(\ref{gptype:A}).
 \item\label{gptype:C} Group generated by the group in (\ref{gptype:A}) and $T$
and not conjugate to a group of type~(\ref{gptype:A}) or (\ref{gptype:B}).
 \item\label{gptype:D} Group generated by the group in (\ref{gptype:C}) and $Q$
and not conjugate to a group of types~(\ref{gptype:A})---(\ref{gptype:C}).
 \item\label{gptype:E} Group of size $108$ generated by $S$, $T$ and $V$.
 \item\label{gptype:F} Group of size $216$ generated by the group in
(\ref{gptype:E}) and an element $P\defby UVU^{-1}$.
 \item\label{gptype:G} Hessian group of size $648$ generated by the group in
(\ref{gptype:E}) and $U$.
 \item\label{gptype:H} Simple group of size $60$ isomorphic to alternating
group \altgroup{5}.
 \item\label{gptype:I} Simple group of size $168$ isomorphic to permutation
group generated by $\left(1234567\right)$, $\left(142\right)\left(356\right)$,
$\left(12\right)\left(35\right)$.
 \item\label{gptype:J} Group of size $180$ generated by the group in
(\ref{gptype:H}) and $W$.
 \item\label{gptype:K} Group of size $504$ generated by the group in
(\ref{gptype:I}) and $W$.
 \item\label{gptype:L} Group $G$ of size $1080$ with its quotient
$G/\defsetspan{W}$ isomorphic to the alternating group \altgroup{6}.
\end{enumerate}
\end{prop}
These groups can be put into the form of the algebraic classification of groups
given in the earlier sections as follows:
\begin{rmk} The list of groups in Proposition~\ref{sl3:prop:classification} can
be subdivided as follows:
\begin{itemize}
 \item Groups of types (\ref{gptype:A}) and (\ref{gptype:B}) are not
irreducible, and so by Theorem~\ref{thm:WExc-trans} cannot induce a weakly-exceptional
singularity.
 \item Groups of types (\ref{gptype:E})---(\ref{gptype:L}) are primitive.
 \item Groups of types (\ref{gptype:C}) and (\ref{gptype:D}) are irreducible monomial.
\end{itemize}
\end{rmk}
\pf Immediate from the lists of generators given above.\qed

Let $G\subset\slgroup{3}{\complex}$ be an irreducible finite subgroup with image
$\bar{G}$ under the natural projection into \pglgroup{3}{\complex}.
\begin{lemma}\label{lemma:sl3:badGroups}
 Assume the singularity the action of $G$ induces is not weakly-exceptional.
Then:
\begin{itemize}
 \item $\bar{G}$ leaves a smooth curve $C\subset\proj{2}$ of degree $2$ invariant.
 \item $G$ is isomorphic to one of \dihgroup{2n},
\altgroup{4}, \symgroup{4}, \altgroup{5} (some $n\geq2$) or to one of their
central extensions by scalar elements.
\end{itemize}
\end{lemma}
\pf By Theorem~\ref{thm:criterion:3D}, $\bar{G}$ must preserve a curve
$C\subset\proj{2}$ of degree $2$. If $C$ is singular, then it must have exactly one isolated
singularity, which is impossible by Lemma~\ref{lemma:singNo}. Therefore, $C$
must be smooth and hence rational, with $\bar{G}$ isomorphic to a finite
irreducible subgroup of \Aut{\proj{1}}.

\pf[ of Theorem~\ref{sl3:thm:main}] Groups of types~(\ref{gptype:A}) and
(\ref{gptype:B}) can be excluded immediately, since they are not irreducible.
Assume the singularity $G$ induces is not weakly-exceptional. Then, comparing
the lists in Proposition~\ref{sl3:prop:classification} and
Lemma~\ref{lemma:sl3:badGroups}, $G$ must be conjugate to one of:
\begin{itemize}
 \item A central extension of
$\left(\cyclgroup{2}\right)^{2}\rtimes\cyclgroup{3}\cong\altgroup{4}$, where both
\cyclgroup{2}-s act diagonally and \cyclgroup{3} permutes
the
basis (such a group is of type~(\ref{gptype:C})).
 \item A central extension of
$\left(\cyclgroup{2}\right)^{2}\rtimes\symgroup{3}\cong\symgroup{4}$, where both
\cyclgroup{2}-s act diagonally and \symgroup{3} permutes the basis
(such a group is of type~(\ref{gptype:D})).
 \item A central extension of \altgroup{5} (such a group is of
type~(\ref{gptype:H}) or (\ref{gptype:J})).
\end{itemize}
First consider the central extensions of \altgroup{4} and \symgroup{4}. Use their presentations from Proposition~\ref{sl3:prop:classification}, and let $\left(x,y,z\right)$ be the corresponding coordinates for $\complex^{3}$. Then the groups have a semi-invariant smooth conic defined by
$x^{2}+y^{2}+z^{2}=0$.
Similarly, it is easy to check (see~\cite[Section~2.9]{Yau-Yu}) that the groups of type~(\ref{gptype:H}) have semi-invariant smooth conics. Since any group of type~(\ref{gptype:J}) is generated by a group of type~(\ref{gptype:H}) and the scalar matrices, these groups will also have semi-invariant smooth conics.

The statement of Theorem~\ref{sl3:thm:main} follows by renaming the remaining
groups to fit with the more widely-used notation.
\qed

\section{Four-dimensional case}\label{sect:4D}
The aim of this section is to prove Theorem~\ref{sl4:thm:main}. By
Theorem~\ref{thm:WExc-trans}, if $G$ is not irreducible, it cannot induce a
weakly-exceptional singularity. Thus, one can restrict attention to the
irreducible groups. Throughout this
section, let $G\subset\slgroup{4}{\complex}$ be a finite irreducible subgroup, and
$\bar{G}$ its image under the natural projection
\deffun{\slgroup{4}{\complex}}{\pglgroup{4}{\complex}}.

In view of Theorem~\ref{thm:criterion:4D}, the only irreducible $\bar{G}$
that are not weakly-exceptional are the finite subgroups of automorphism
groups
of surfaces of degree $2$ or $3$ and \altgroup{5} with one specific action.
When considering automorphisms of a surface $S$, without loss of generality one
can assume that there is no
$\bar{G}$-invariant surface $S'$ of smaller degree.
\begin{lemma} Let $G\subset\slgroup{4}{\complex}$ be a finite irreducible group, $\bar{G}\subset\pglgroup{4}{\complex}$ its projection, and let $S\subset\proj{3}$ be a $\bar{G}$-invariant surface of degree
minimal among the degrees of all $\bar{G}$-invariant surfaces. Then either
$\deg{S}\geq4$ or $S$ is smooth.
\end{lemma}
\pf Since $G$ is irreducible, $\deg{S}\geq2$. If $\deg{S}=2$ and $S$ is singular,
then either $S$ has exactly $1$ isolated singularity, or $S$ is a union of two
planes and thus has a singular line, which must then be $\bar{G}$-invariant.
Both cases are impossible (by Lemma~\ref{lemma:singNo} and the irreducibility of $G$), so
if $\deg{S}=2$ then $S$ must be smooth.

If $\deg{S}=3$, and $S$ is not irreducible, then either it is the union of a
plane and an irreducible quadric surface (each of which must thus be a
$\bar{G}$-invariant surface of smaller degree, contradicting the hypothesis) or
$S$ is the union of $3$ distinct planes, whose intersection gives a point fixed
by all of $\bar{G}$ (stopping $G$ from being irreducible). Hence $S$ is
irreducible.

Assume $S$ has non-isolated singularities, with $C$ being the union of all
singular curves on $S$. Then, one can easily see that $C$ is a line. Since
$\bar{G}(S)=S$, $\bar{G}(C)=C$, and so there exists a $\bar{G}$-invariant line,
contradicting irreducibility of $G$. Therefore if $\deg{S}=3$ then $S$ must have
at worst isolated singularities.

If $\deg{S}=3$ and $S$ is singular with only isolated singularities,
then by~\cite{Bruce-Wall}, the singularity types form one of the following
collections:
$\left(A_{1}\right)$,
$\left(2A_{1}\right)$,
$\left(A_{1},A_{2}\right)$,
$\left(3A_{1}\right)$,
$\left(A_{1},A_{3}\right)$,
$\left(2A_{1},A_{2}\right)$,
$\left(4A_{1}\right)$,
$\left(A_{1},A_{4}\right)$,
$\left(2A_{1},A_{3}\right)$,
$\left(A_{1},2A_{2}\right)$,
$\left(A_{1},A_{5}\right)$.
Since by Lemma~\ref{lemma:singNo}, $S$ cannot have exactly $1$, $2$ or $3$
singularities of any given type, $S$ has to have $4$~$A_{1}$ singularities.
Since there is only one such surface (see, for example, \cite{Bruce-Wall}), $S$
must be the Cayley cubic, defined (in some basis) by
\[S=\defset{\left(x:y:u:v\right)\in\proj{3}}{xyu+xyv+xuv+yuv=0}\]

But in this case $S$ contains exactly $9$ lines, $6$ of which are going through
pairs of singular points and the other $3$ defined by 
\[x+y=0=u+v,\ x+u=0=y+v\ \mbox{and}\ x+v=0=y+u\]
These last three lines are coplanar and must be mapped to each other by all of
$\bar{G}$. Therefore, $\bar{G}$ preserves a plane, contradicting the
irreducibility assumption for $G$. Thus if $S$ is a cubic surface, then it must be
smooth.
\qed

Summarising, if $\bar{G}$ is irreducible but the singularity it
induces is not weakly-exceptional, any $\bar{G}$-invariant surface $S$ of degree at most $3$ must be
a smooth surface of degree $2$ or $3$. These cases will be considered separately
in the next two sections.

\pf[ of Theorem~\ref{sl4:thm:main}] By Theorem~\ref{thm:criterion:4D} and
the discussion above, the theorem is an immediate consequence of
Lemmas~\ref{lemma:sl4:cubic:smooth} and~\ref{lemma:sl4:quadric}.\qed

\subsection{If $S$ is a smooth $\bar{G}$-invariant cubic surface}
This section is devoted to proving the following lemma:
\begin{lemma}\label{lemma:sl4:cubic:smooth}
If $G\subset\slgroup{4}{\complex}$ is a finite irreducible subgroup, and
$\bar{G}$ its projection to $\pglgroup{4}{\complex}$. Also assume that
there is no $\bar{G}$-invariant quadric surface, and
$S\subset\proj{3}$ is a smooth $\bar{G}$-invariant cubic surface. Then $G$ must
be isomorphic to a central extension of one of:
\begin{itemize}
 \item $\left(\left(\cyclgroup{3}\right)^{3}
\rtimes\cyclgroup{2}\right) \rtimes\cyclgroup{2}$. This produces a
monomial action.
 \item $\left(\cyclgroup{3}\right)^{3}\rtimes\symgroup{4}$. This
produces a monomial action.
\end{itemize}
by scalar elements, acting as described below. Comparisons with both of these
isomorphism classes are indeed necessary.
\end{lemma}
As stated before, $\bar{G}\subset\Aut{S}$ is a finite subgroup, so by \cite{Hosoh},
$\bar{G}$
must belong to one of the following isomorphism classes:
\begin{enumerate}
 \item\label{groups:cubic:cyclic} \defsetshort{e}, \cyclgroup{2},
\cyclgroup{4}, \cyclgroup{8}.
 \item $\left(\cyclgroup{2}\right)^{2}$, \symgroup{3}, $\symgroup{3}\times\cyclgroup{2}$.
 \item \symgroup{4}.
 \item $\left(\cyclgroup{3}\right)^{2}\rtimes\cyclgroup{2}$.
 \item $\left(\left(\cyclgroup{3}\right)^{3}
\rtimes\cyclgroup{2}\right)
\rtimes\cyclgroup{2}$.
 \item \symgroup{5}.
 \item $\left(\cyclgroup{3}\right)^{3}\rtimes\symgroup{4}$.
\end{enumerate}

\subsubsection{Case $\bar{G}$ cyclic}
The groups in (\ref{groups:cubic:cyclic}) are all cyclic, so their extensions by
scalar elements do not act irreducibly.
\subsubsection{Case $\bar{G}$ dihedral}
This case covers the isomorphism classes $\left(\cyclgroup{2}\right)^{2}$,
\symgroup{3} and $\symgroup{3}\times\cyclgroup{2}$. The dihedral groups and
their extensions by scalar elements do not have any irreducible $4$-dimensional
representations, so these groups cannot act irreducibly.
\subsubsection{Case $\bar{G}\cong\symgroup{4}$} This group by itself has no
$4$-dimensional irreducible representations, while its central extension has
(up to a choice of a root of unity) only one such, which preserves a quadric
surface (see the twisted diagonal actions in Lemma~\ref{lemma:sl4:quadric}).
\subsubsection{Cases $\bar{G}\cong\left(\cyclgroup{3}\right)^{2}
\rtimes\cyclgroup{2}$ or
$\left(\left(\cyclgroup{3}\right)^{3}\rtimes\cyclgroup{2}\right)
\rtimes\cyclgroup{2}$}
For convenience, write $\bar{G}'=\left(\cyclgroup{3}\right)^{2}
\rtimes\cyclgroup{2}$ and
$\bar{G}''=\left(\left(\cyclgroup{3}\right)^{3}\rtimes\cyclgroup{2}\right)
\rtimes\cyclgroup{2}$, with all the notation following in the obvious manner
(i.e.\ write $G'$ for the lift of $\bar{G}'$ to \slgroup{4}{\complex}, etc.).
Using the notation from~\cite{Hosoh}, these two cases correspond to groups
$\bar{G}'=G_{54}^{9}/C\left(G_{54}^{9}\right)$ and
$\bar{G}''=G_{54}^{9}\rtimes\cyclgroup{2}$. This means there exist non-trivial
elements
$\bar{\alpha},\bar{\beta},\bar{\gamma},\bar{\delta},\bar{\epsilon}
\in\bar{G}''$,
such that $\bar{\alpha},\bar{\beta},\bar{\gamma},\bar{\delta}$ generate
$G_{54}^{9}$, with $\bar{\alpha}$ generating its centre,
$\bar{\alpha}^{3}=\bar{\beta}^{3}=\bar{\gamma}^{3}=\bar{\delta}^{2}=\bar{\epsilon}^{2}=\mbox{id}$. Let
$\alpha,\beta,\gamma,\delta,\epsilon$ be lifts of
$\bar{\alpha},\bar{\beta},\bar{\gamma},\bar{\delta}$ (respectively) to
\slgroup{4}{\complex}.

Let $h_{1}\defby\alpha^{3}$, $h_{2}\defby\beta^{3}$. $\bar{\alpha},\bar{\beta}$
commute, so say $\beta\alpha=\alpha\beta h_{3}$. By the structure of the lift,
$h_{i}$ are scalar matrices of order $1$, $2$ or $4$. Then
\[h_{1}^{3}h_{2}=\left(\alpha^{2}\beta\alpha\right)^{3}=
\left(\beta h_{1}h_{3}\right)^{3}=h_{2}h_{1}^{3}h_{3}^{3}\]
and so $h_{3}=\mbox{id}$. Similarly, get $\alpha,\beta,\gamma$ all commuting.
Hence the corresponding matrices can all be taken to be diagonal (by choosing a
suitable basis). It is then easy to see that $\delta$ and $\epsilon$ must act
as elements of a central extension of \symgroup{4} permuting the basis.

Since $\bar{G}''$ has only one normal subgroup of index $2$, and $\bar{G}''$ has
no centre (otherwise $\bar{G}''/C\left(\bar{G}''\right)$ would be on the list of
groups acting on a cubic surface),
$\bar{\delta}\bar{\epsilon}\neq\bar{\epsilon}\bar{\delta}$. Therefore, up to
conjugation, $\delta$ interchanges the first and the second basis vectors, and
$\epsilon$ interchanges the first basis vector with the third one and the second
basis vector with the fourth one.

This means that $G'$ is not irreducible, while $G''$ is irreducible and (up
to conjugation) is generated by
\[
\left(\begin{array}{cccc}
\zeta_{3}&0&0&0\\ 0&1&0&0\\ 0&0&1&0\\ 0&0&0&\zeta_{3}^{-1}
\end{array}\right),\ 
\left(
\begin{array}{cccc}
1&0&0&0\\ 0&\zeta_{3}&0&0\\ 0&0&1&0\\ 0&0&0&\zeta_{3}^{-1}
\end{array}
\right),\ 
\left(\begin{array}{cccc}
1&0&0&0\\ 0&1&0&0\\ 0&0&\zeta_{3}&0\\ 0&0&0&\zeta_{3}^{-1}
\end{array}\right),\ 
\]
\[
\zeta_{8}\left(\begin{array}{cccc}
0&1&0&0\\ 1&0&0&0\\ 0&0&1&0\\ 0&0&0&1
\end{array}\right),\ 
\left(\begin{array}{cccc}
0&0&1&0\\ 0&0&0&1\\ 1&0&0&0\\ 0&1&0&0
\end{array}\right)
\]

This group leaves (for example) the cubic polynomial $x^{3}+y^{3}+z^{3}+w^{3}$
(in coordinates $\left(x,y,z,w\right)$ for $\complex^{4}$) semi-invariant, and
by direct computation, one sees that the group does not have a semi-invariant
quadric surface. This action of $G''$ is monomial.
\subsubsection{Case $\bar{G}\cong\symgroup{5}$}
According to~\cite[\S100]{Segre}, this is the automorphism group of the
irreducible diagonal cubic surface
\[
S=\defset{\left(x_{0}:x_{1}:x_{2}:x_{3}:x_{4}\right)\in\proj{4}}{
\begin{array}{c} x_{0}^{3}+x_{1}^{3}+x_{2}^{3}+x_{3}^{3}+x_{4}^{3}=0,\\
x_{0}+x_{1}+x_{2}+x_{3}+x_{4}=0\end{array}}
\]
which immediately implies that there exists a $\bar{G}$-invariant quadric
surface. For the group's action on it, see Lemma~\ref{lemma:sl4:quadric}.
\subsubsection{Case
$\bar{G}\cong\left(\cyclgroup{3}\right)^{3}\rtimes\symgroup{4}$}
As stated in \cite[\S100]{Segre}, this group acts by permuting the basis vectors of
$\complex^{4}$ arbitrarily and multiplying them by arbitrary cube roots of
unity. Hence (up to conjugation) $G$ is a central extension of such a group by
scalar elements.

This group clearly leaves the cubic polynomial $x^{3}+y^{3}+z^{3}+w^{3}$ (in
coordinates $\left(x,y,z,w\right)$) semi-invariant, and by direct computation,
one sees that the group does not have a semi-invariant quadric surface. The
action of this group is monomial.

\subsection{If $S$ is a smooth $\bar{G}$-invariant quadric surface}\label{sect:sl4:quadric}
Let the group $\bar{G}\subset\pglgroup{4}{\complex}$ be a finite irreducible
subgroup, $G$ its lift to \slgroup{4}{\complex} and $S\subset\proj{3}$ a smooth
$\bar{G}$-invariant quadric surface. This section is devoted to compiling a
list of the possible values that $G$ (equivalently, $\bar{G}$) can take in this
situation. The final list is presented in Lemma~\ref{lemma:sl4:quadric}.

In this case there exists a basis for \proj{3}, in which $S$ is the image of
the Segre embedding
of $\proj{1}\times\proj{1}$ into \proj{3}. The image of this embedding is
\[S=\defset{\left(x:y:u:v\right)\in\proj{3}}{xv-yu=0}\]
This implies that the subgroup of \pglgroup{4}{\complex} preserving $S$ is
isomorphic to $\left(\pglgroup{2}{\complex}\right)^{2}\rtimes\cyclgroup{2}$,
with the action that can be described in the notation given in
Sections~\ref{notation} and~\ref{sect:quadric-action}. This notation will be used throughout the rest of the paper.

Now one needs to put the action of any finite subgroup of \Aut{S} into one of the four categories:
\begin{itemize}
 \item Not irreducible.
 \item Irreducible monomial.
 \item Irreducible non-monomial imprimitive.
 \item Primitive (hence irreducible).
\end{itemize}
This can be done directly by looking at the representations used to build the
action. To make the explanations more simple, it will be assumed that $H_{1}$
and $H_{2}$ (see Section~\ref{sect:quadric-action}) both contain a non-scalar diagonal matrix. This will mean that any
proper invariant subspace must have a basis, which is a subset of the chosen
basis for $\complex^{4}$. It is easy to check that for all the actions used in this section, there exists a basis for $\complex^{4}$ in which the action contains such a matrix. 
Now fix this basis for the remainder of this section.

\subsubsection{Irreducibility:}
Assume first that $\bar{G}$ does not interchange the rulings. Then it can be seen that in order for $G$ to be
irreducible,
both $H_{1}$ and $H_{2}$ need to be irreducible: if $H_{1}$ is not irreducible,
then (in the matrix presentation of \proj{3}) the ``rows'' of \proj{3} produce
invariant lines, i.e.
\[\defsetshort{x=y=0},\defsetshort{u=v=0}\subset
\defsetshort{\left(x:y:u:v\right)\in\proj{3}}\]
are invariant subspaces. Similarly, if $H_{2}$ is not irreducible, then there
exist $H$-invariant lines corresponding to the ``columns'' of \proj{3}, i.e.
\[\defsetshort{x=u=0},\defsetshort{y=v=0}\subset
\defsetshort{\left(x:y:u:v\right)\in\proj{3}}\]
are invariant subspaces. Even when both $H_{i}$ are irreducible, $G$ can still
fail to be irreducible, but this can only occur if both $H_{i}$ are dihedral
groups, with their order $2$ elements coupled (see the end
of Section~\ref{sect:quadric-action}).

This means that for the action of $G$ to be irreducible, need (for any $i\neq
j\in\defsetshort{1,2})$) one of the following to hold:
\begin{itemize}
 \item $H_{i},H_{j}\in\defsetshort{\altgroup{4},\symgroup{4},\altgroup{5}}$
 \item $H_{i}\in\defsetshort{\altgroup{4},\symgroup{4},\altgroup{5}}$,
$H_{j}=\dihgroup{2n}$
 \item $H_{1}=\dihgroup{4m}$, $H_{2}=\dihgroup{4n}$ and the action contains
elements of the
form $\quadricaction{a^{l}}{b}$ or $\quadricaction{b}{a^{k}}$, not just $\quadricaction{a^{l}b}{a^{k}b}$ ($k,l\in\integers$).
\end{itemize}

If $\bar{G}$ does interchange the rulings, then assume first that $\tau\in\bar{G}$ (where $\tau$ is the involution interchanging the rulings, see Section~\ref{sect:quadric-action}). Then $\tau$ must keep one of the diagonals of the
matrix form of \proj{3} invariant, i.e.
\[\defsetshort{x=v=0},\defsetshort{y=u=0}\subset
\defsetshort{\left(x:y:u:v\right)\in\proj{3}}\]
are left invariant by $\tau$. However, if $\tau\in\bar{G}$, then $H_{1}\cong
H_{2}$, so the only two cases where it can potentially make a previously not irreducible
action into an irreducible one are
when $H_{1}\cong H_{2}\cong\dihgroup{2n}$ (the action is
not irreducible, as $\quadricaction{a^{l}b}{a^{k}b}$ leaves the same two proper subspaces
invariant as $\tau$ does) or when $H_{1}\cong H_{2}\cong\cyclgroup{n}$
(the action is not irreducible, as $\tau$ only permutes pairs of coordinates,
so the group has a pair of distinct invariant proper subspaces)

If $\bar{G}$ interchanges the rulings, but $\tau\not\in\bar{G}$, then $\sigma\tau\in\bar{G}$ for some automorphism $\sigma$ that preserves the ruling, with $\sigma\not\in H$, $\sigma^{2}\in H$. Consider $\hat{G}$ generated by $H$, $\tau$ and $\sigma$, with $\hat{K_{i}},\hat{H_{i}},\hat{H}$ defined similarly to $K_{i},H_{i},H$. Clearly, $\bar{G}\subset\hat{G}$, $H\subset\hat{H}$, etc. 
\begin{itemize}
 \item If $H_{i}\cong\cyclgroup{m}$, then either $\hat{H_{i}}\cong\cyclgroup{2m}$ or $\hat{H_{i}}\cong\dihgroup{2m}$ with $\sigma=\quadricaction{a^{l}b}{a^{k}b}\in\hat{H}$. (some $k,l\in\integers$). In either case this makes $\hat{G}$ not irreducible, and hence $\bar{G}$ not irreducible (by the discussion above).
 \item If $H_{i}\cong\dihgroup{2m}$, then it is easy to see that $\hat{H_{i}}\cong\dihgroup{2n}$ (for some $n>m$). Similarly to the considerations above, this only influences the irreducibility of $\bar{G}$ if $\sigma$ can be taken to be of the form $\quadricaction{a^{l}b}{a^{k}}$, making the group action irreducible.
\end{itemize}

The discussion above has been summarised in
Table~\ref{sl4:quadrics:transitivity}.
\begin{table}[htbp]
\caption{Irreducibility of $G$}\label{sl4:quadrics:transitivity}
\begin{center}
\begin{tabular}{cr|ccc}
&&\multicolumn{3}{c}{$H_{1}$}\\
&&\altgroup{4}, \symgroup{4}, \altgroup{5}& \dihgroup{2m}&
\cyclgroup{m}\\
\hline
\multirow{3}{*}{$H_{2}$}&\altgroup{4}, \symgroup{4},
\altgroup{5}&Irreducible&Irreducible&Not irreducible\\
&\dihgroup{2m}&Irreducible&Depends on action&Not irreducible\\
&\cyclgroup{n}&Not irreducible&Not irreducible&Not irreducible
\end{tabular}
\end{center}
\end{table}

\subsubsection{Primitivity:}
Assume the action of $G$ is irreducible. Again, the notation from Section~\ref{sect:quadric-action} will be used throughout. Then the question of the action being monomial, imprimitive non-monomial or primitive becomes of interest.

By direct computation, it is easy to see that in most cases, the place of $G$ in this
classification depends on the matrices in $H_{i}$ that have $3$ or more
non-zero entries, and how these matrices are combined in $G$, i.e.\ on the
isomorphism \deffunname{\alpha}{H_{1}/\hat{K}_{1}}{H_{2}/\hat{K}_{2}} (as
defined above). The only exception occurs when $H_{i}\cong\dihgroup{2m}$, $G$ interchanges the ruling, but $\tau\not\in\bar{G}$ --- in this case, the automorphism $\sigma$ needs to be considered.

With this in mind, direct computation provides the following criteria (putting
$i\neq j\in\defsetshort{1,2}$):
\begin{itemize}
 \item If $H_{1},H_{2}$ dihedral and $G$ irreducible, then $G$ acts monomially.
 \item If $H_{i}\in\defsetshort{\altgroup{4}, \symgroup{4}, \altgroup{5}}$ and
$H_{j}=\dihgroup{2n}$, then the action of $G$ is non-monomial imprimitive.
 \item If $H_{1}\cong H_{2}\cong\altgroup{5}$, then the action of $G$ is
primitive.
 \item If $H_{1},H_{2}\in\defsetshort{\altgroup{4},\symgroup{4}}$ and the
$3$-cycles are not coupled, then the action of $G$ is primitive.
 \item If $H_{1}\cong H_{2}\cong\altgroup{4}$ and the $3$-cycles are coupled,
then the action of $G$ is monomial.
 \item If $H_{i}\cong\symgroup{4}$, $H_{j}\cong\altgroup{4}$ and the
$3$-cycles are coupled, then the action is imprimitive non-monomial.
 \item If $H_{i}\cong\symgroup{4}$, $H_{j}\cong\symgroup{4}$, the $3$-cycles
are coupled and the odd permutations are coupled, then the action is monomial.
 \item If $H_{i}\cong\symgroup{4}$, $H_{j}\cong\symgroup{4}$, the $3$-cycles
are coupled, but the odd permutations are not coupled, then the action is
primitive.
\end{itemize}
This list is clearly not exhaustive, but it is sufficient for determining the
nature of all the groups below.

\subsubsection{Possible isomorphism classes of $\bar{G}$}
Since $\bar{G}$ is a finite group leaving a smooth quadric $S$ invariant, its
action must be equal (as shown in Section~\ref{prelim:invars}) to a suitable of
one
of the finite automorphism groups of a smooth $2$-dimensional quadric. Thus
$\bar{G}$ must be conjugate to the image of one of the finite groups given
in~\cite[Theorem~4.9]{D-I}.

In order to make the structure of each of the groups slightly more explicit,
the group structure will also be given in the notation
\[\left(H_{1},K_{1},H_{2},K_{2}\right)_{\alpha},\] where $H_{i}$, $K_{i}$ are as
before, and $\alpha$ is the gluing isomorphism between $H_{1}/\hat{K}_{1}$ and
$H_{2}/\hat{K}_{2}$. Where only one such isomorphism exists, $\alpha$ will be
omitted. For each isomorphism class, several representations of the group can
be chosen. However, 
it is clear that the different faithful representations will differ by at most an outer
automorphism, so all the properties that are of interest in this discussion
will be the same for all of them. Therefore, for each isomorphism class, any
faithful representation of $H_{i}$ can be chosen. For any group
$\left(H_{1},K_{1},H_{2},K_{2}\right)_{\alpha}$, there also exists a group
$\left(H_{2},K_{2},H_{1},K_{1}\right)_{\alpha^{-1}}$, which corresponds to the
same group with the components of the ruling of the quadric swapped. These two
groups are conjugate to each other.

\begin{lemma}\label{lemma:sl4:quadric}
If $\bar{G}\subset\pglgroup{4}{\complex}$ is a finite irreducible subgroup,
$G$ its lift to \slgroup{4}{\complex} and $S\subset\proj{3}$ a smooth
$\bar{G}$-invariant quadric surface, then $\bar{G}$ must be
conjugate to one of the groups in the list below. From the way the groups are
constructed, it is easy to see that all the groups do need to be in the list.
\end{lemma}

First assume $\bar{G}$ leaves the ruling of
$\proj{1}\times\proj{1}$ invariant. Bearing the analysis above in mind, for
$G$ to be irreducible, it must be conjugate to one of the following:
\begin{enumerate}
 \item Product subgroups $\left(H_{1},H_{1},H_{2},H_{2}\right)\cong H_{1}\times
H_{2}$ for some finite groups $H_{i}\in\Aut{\proj{1}}$. Taking different
choices for $H_{1},H_{2}$,
get the following groups of the form $H_{1}\times H_{2}$:
\begin{enumerate}
\item 9 primitive groups when
$H_{1},H_{2}\in\defsetshort{\altgroup{4},\symgroup{4},\altgroup{5}}$.
\item 3 families of non-monomial imprimitive groups $\dihgroup{2m}\times H_{2}$,
where $H_{2}\in\defsetshort{\altgroup{4},\symgroup{4},\altgroup{5}}$.
\item 1 family of monomial groups $\dihgroup{2m}\times\dihgroup{2n}$.
\end{enumerate}
 \item Twisted diagonal subgroups $\left(H_{1},1,H_{1},1\right)_{\alpha}\cong
H_{1}$ for some finite group $H_{1}\in\Aut{\proj{1}}$. This gives $3$ families
of groups,
indexed by the choice of isomorphism $\alpha$. They are:
\begin{enumerate}
 \item Monomial groups isomorphic to \altgroup{4} or \symgroup{4}.
 \item Primitive groups isomorphic to \altgroup{5}.
\end{enumerate}
The twisted diagonal groups isomorphic to the dihedral groups do not act
irreducibly, as the relevant central extensions do not have any
$4$-dimensional irreducible representations.

 \item\label{usedin:prim}
$\frac{1}{2}\left[\symgroup{4}\times\symgroup{4}\right]
\cong\left(\altgroup{4}\times\altgroup{4}\right)\rtimes\cyclgroup{2}
\cong\left(\symgroup{4},\altgroup{4},\symgroup{4},\altgroup{4}\right)$,
a primitive group generated by elements corresponding to
$\quadricaction{[12][34]}{\mbox{id}}$,
$\quadricaction{\mbox{id}}{[12][34]}$,
$\quadricaction{[123]}{\mbox{id}}$,
$\quadricaction{\mbox{id}}{[123]}$ and
$\quadricaction{(12)}{(12)}$.

 \item\label{usedin:imprim:nonmono}
$\frac{1}{2}\left[\dihgroup{2m}\times\symgroup{4}\right]
\cong\left(\cyclgroup{m}\times\altgroup{4}\right)\rtimes\cyclgroup{2}
\cong\left(\dihgroup{2m},\cyclgroup{m},\symgroup{4},\altgroup{4}\right)$,
a family of imprimitive non-monomial groups generated by
$\quadricaction{a_{m}}{\mbox{id}}$,
$\quadricaction{\mbox{id}}{[12][34]}$,
$\quadricaction{\mbox{id}}{[123]}$ and
$\quadricaction{b}{(12)}$.

 \item
$\frac{1}{2}\left[\dihgroup{4m}\times\symgroup{4}\right]
\cong\left(\dihgroup{2m}\times\altgroup{4}\right)\rtimes\cyclgroup{2}
\cong\left(\dihgroup{4m},\dihgroup{2m},\symgroup{4},\altgroup{4}\right)$
$\left(m\geq2\right)$,
a family of imprimitive non-monomial groups generated by the action of
$\quadricaction{a_{2m}^{2}}{\mbox{id}}$,
$\quadricaction{b}{\mbox{id}}$,
$\quadricaction{\mbox{id}}{[12][34]}$,
$\quadricaction{\mbox{id}}{[123]}$ and
$\quadricaction{a_{2m}^{m}}{(12)}$.

 \item
$\frac{1}{6}\left[\dihgroup{6m}\times\symgroup{4}\right]
\cong\left(\cyclgroup{m}\times \mathbb{V}_{4}\right)\rtimes\symgroup{3}
\cong\left(\dihgroup{6m},\cyclgroup{m},\symgroup{4},\mathbb{V}_{4}\right)$,
a family of imprimitive non-monomial groups generated by
$\quadricaction{a_{3m}^{3}}{\mbox{id}}$,
$\quadricaction{\mbox{id}}{[12][34]}$,
$\quadricaction{\mbox{id}}{(13)(24)}$,
$\quadricaction{a_{3m}^{m}}{[123]}$ and
$\quadricaction{b}{(12)}$.

 \item\label{usedin:imprim}
$\frac{1}{6}\left[\symgroup{4}\times\symgroup{4}\right]
\cong \left(\mathbb{V}_{4}\times \mathbb{V}_{4}\right)\rtimes\symgroup{3}
\cong\left(\symgroup{4},\mathbb{V}_{4},\symgroup{4},\mathbb{V}_{4}\right)$,
a monomial group generated by
$\quadricaction{[12][34]}{\mbox{id}}$,
$\quadricaction{(13)(24)}{\mbox{id}}$,
$\quadricaction{\mbox{id}}{[12][34]}$,
$\quadricaction{\mbox{id}}{(13)(24)}$,
$\quadricaction{[123]}{[123]}$ and
$\quadricaction{(12)}{(12)}$.

 \item
$\frac{1}{3}\left[\altgroup{4}\times\altgroup{4}\right]
\cong\left(\mathbb{V}_{4}\times \mathbb{V}_{4}\right)\rtimes\cyclgroup{3}
\cong\left(\altgroup{4},\mathbb{V}_{4},\altgroup{4},\mathbb{V}_{4}\right)$,
a monomial group generated by
$\quadricaction{[12][34]}{\mbox{id}}$,
$\quadricaction{(13)(24)}{\mbox{id}}$,
$\quadricaction{\mbox{id}}{[12][34]}$,
$\quadricaction{\mbox{id}}{(13)(24)}$ and
$\quadricaction{[123]}{[123]}$.

 \item
$\frac{1}{2}\left[\dihgroup{2m}\times\dihgroup{4n}\right]
\cong\left(\cyclgroup{m}\times\dihgroup{2n}\right)\rtimes\cyclgroup{2}
\cong\left(\dihgroup{2m},\cyclgroup{m},\dihgroup{4n},\dihgroup{2n}\right)$
(for $m,n\geq2$),
a family of monomial groups generated by 
$\quadricaction{a_{m}}{\mbox{id}}$,
$\quadricaction{\mbox{id}}{a_{2n}^{2}}$,
$\quadricaction{\mbox{id}}{b}$ and
$\quadricaction{b}{a_{2n}^{n}}$.

 \item\label{usedin:trans}
$\frac{1}{4}\left[\dihgroup{4m}\times\dihgroup{4n}\right]_{\alpha}
\cong\left(\cyclgroup{m}\times \cyclgroup{n}\right)\rtimes\dihgroup{4}
\cong\left(\dihgroup{4m},\cyclgroup{m},\dihgroup{4n},\cyclgroup{n}\right)_{
\alpha }$ (where
$\alpha(b)=a_{2n}^{n}$, $\alpha(a_{2m}^{m})=b$),
a family of monomial groups generated by
$\quadricaction{a_{2m}^{2}}{\mbox{id}}$,
$\quadricaction{\mbox{id}}{a_{2n}^{2}}$,
$\quadricaction{a_{2m}^{m}}{b}$ and
$\quadricaction{b}{a_{2n}^{n}}$.

 \item
$\frac{1}{2}\left[\dihgroup{4m}\times\dihgroup{4n}\right]
\cong\left(\dihgroup{2m}\times\dihgroup{2n}\right)\rtimes\cyclgroup{2}
\cong\left(\dihgroup{4m},\dihgroup{2m},\dihgroup{4n},\dihgroup{2n}\right)$
(for $m,n\geq2$),
a family of monomial groups generated by 
$\quadricaction{a_{2m}^{2}}{\mbox{id}}$,
$\quadricaction{b}{\mbox{id}}$,
$\quadricaction{\mbox{id}}{a_{2n}^{2}}$,
$\quadricaction{\mbox{id}}{b}$ and
$\quadricaction{a_{2m}^{m}}{a_{2n}^{n}}$.

\suspend{enumerate}

Now assume that $\bar{G}$ does interchange the rulings (via an element $\sigma\circ\tau$). Then the normal subgroup of $\bar{G}$
fixing the ruling of the quadric must also be a group of automorphisms of $S$.
Furthermore, it must have $H_{1}\cong H_{2}$ and $K_{1}\cong K_{2}$, since
conjugation by $\sigma\tau$ provides the two isomorphisms. That means that $\bar{G}$
can be isomorphic to one of the following groups:
\resume{enumerate}

 \item $\left(H_{1}\times H_{1}\right)\rtimes2\cong
\left(H_{1},H_{1},H_{1},H_{1}\right)\rtimes\cyclgroup{2}$
($H_{1}\in\Aut{\proj{1}}$).
Taking different choices for $H_{1}$ and bearing in mind that choosing
$H_{1}$ to be \cyclgroup{n} produces a group that is not irreducible (see discussion
above), get $2$ families of monomial groups
\begin{enumerate}
 \item
$\left(\dihgroup{2n}\times\dihgroup{2n}\right)\rtimes\cyclgroup{2}$
\suspend{enumerate}
and $3$ families of primitive groups, all of them indexed by the possible
involutions acting on $H_{1}$:
\resume{enumerate}
 \item $\left(\altgroup{4}\times\altgroup{4}\right)\rtimes\cyclgroup{2}$.
 \item $\left(\symgroup{4}\times\symgroup{4}\right)\rtimes\cyclgroup{2}$.
 \item $\left(\altgroup{5}\times\altgroup{5}\right)\rtimes\cyclgroup{2}$.
\end{enumerate}

 \item $H_{1}\rtimes\cyclgroup{2}\cong \left(H_{1},1,H_{1},1\right)_{\alpha}$
($H_{1}\in\Aut{\proj{1}}$). This gives $3$ families of groups,
indexed by the choice of isomorphism $\alpha$. They are:
\begin{enumerate}
 \item Monomial groups isomorphic to $\dihgroup{4n}\rtimes\cyclgroup{2}$.
 \item Monomial groups isomorphic to $\altgroup{4}\rtimes\cyclgroup{2}$ or
$\symgroup{4}\rtimes\cyclgroup{2}$.
 \item Primitive groups isomorphic to $\altgroup{5}\rtimes\cyclgroup{2}$.
\end{enumerate}

 \item
$\left(\left(\altgroup{4}\times\altgroup{4}
\right)\rtimes\cyclgroup{2}\right)\rtimes\cyclgroup{2}
\cong\left(\symgroup{4},\altgroup{4},\symgroup{4},\altgroup{4}
\right)\rtimes\cyclgroup{2}$,
a family of primitive groups.

 \item
$\left(\left(\mathbb{V}_{4}\times
\mathbb{V}_{4}\right)\rtimes\symgroup{3}\right)\rtimes\cyclgroup{2}
\cong\left(\symgroup{4},\mathbb{V}_{4},\symgroup{4},\mathbb{V}_{4}
\right)\rtimes\cyclgroup{2}$,
a family of monomial groups.

 \item\label{usedin:imprim:mono}
$\left(\left(\mathbb{V}_{4}\times
\mathbb{V}_{4}\right)\rtimes\cyclgroup{3}\right)\rtimes\cyclgroup{2}
\cong\left(\altgroup{4},\mathbb{V}_{4},\altgroup{4},\mathbb{V}_{4}
\right)\rtimes\cyclgroup{2}$,
a family of monomial groups.

 \item
$\left(\left(\cyclgroup{m}\times\cyclgroup{m}\right)\rtimes\dihgroup{4}
\right)\rtimes\cyclgroup{2}
\cong\left(\dihgroup{4m},\cyclgroup{m},\dihgroup{4m},
\cyclgroup{m}\right)_{\alpha}\rtimes\cyclgroup{2}$
(where
$\alpha(b)=a_{2m}^{m}$, $\alpha(a_{2m}^{m})=b$),
a family of monomial groups.

 \item
$\left(\left(\dihgroup{2m}\times\dihgroup{2m}\right)\rtimes\cyclgroup{2}
\right)\rtimes\cyclgroup{2}
\cong\left(\dihgroup{4m},\dihgroup{2m},\dihgroup{4m},\dihgroup{2m}
\right)\rtimes\cyclgroup{2}$
($m\geq2$),
a family of monomial groups.
\end{enumerate}

This concludes the proof of Theorem~\ref{sl4:thm:main}, with the original list
of groups being a re-ordering of this list combined with the list in
Lemma~\ref{lemma:sl4:cubic:smooth}.

As a final remark, it is worth noting that to make the proof more transparent
while keeping the size of it reasonable, the presentations for some of the isomorphism classes were chosen to match the structure of the proof more closely, instead of being the more common ones, for example:
\begin{eg}
When one thinks about finite groups in \slgroup{4}{\complex} that leave a smooth
quadric invariant, one of the first examples that come to mind is the monomial
group $G=\left(\cyclgroup{2}\right)^{3}\rtimes\symgroup{4}$, where the normal
subgroup $H=\left(\cyclgroup{2}\right)^{3}$ acts
diagonally, and the symmetric group permutes the basis. Here, one of the
size two subgroups of $H$ acts by scalar matrices, so the action on the
projective quadric surface is isomorphic to
$\left(\cyclgroup{2}\right)^{2}\rtimes\symgroup{4}\cong\left(\left(\mathbb{V}_{4
}\times
\mathbb{V}_{4}\right)\rtimes\cyclgroup{3}\right)\rtimes\cyclgroup{2}$, which is
in
position~(\ref{usedin:imprim:mono}) in the list above.
\end{eg}

\def\cprime{$'$}

\end{document}